\documentclass[11pt,reqno]{article}
\usepackage{amsmath, amssymb, amsthm}

\numberwithin{equation}{section}

\newtheorem{theorem}{Theorem}[section]

\newtheorem{proposition}[theorem]{Proposition}
\newtheorem{corollary}[theorem]{Corollary}
\newtheorem{lemma}[theorem]{Lemma}

\newtheorem{fact}[theorem]{Fact}
\newtheorem*{general Gromov'}{Corollary \ref{general Gromov}$'$}

\def \proof {\noindent {\bf Proof.}\ \ }
\def \remark {\noindent {\bf Remark.}\ \ }

\def \endproof {{\mbox{}\nolinebreak\hfill\rule{2mm}{2mm}\par\medbreak}}

\def \R {\mathbb{R}}

\def \P {\mathbb{P}}

\def \NN {\mathcal{N}}

\def \OO {\mathcal{O}}

\def \a {\alpha}
\def \b {\beta}
\def \g {\gamma}
\def \e {\varepsilon}

\def \d {\delta}

\def \s {\sigma}

\def \< {\langle}
\def \> {\rangle}

\def \rank {{\rm rank }}

\def \Lip {{\rm Lip}}

\begin{document}
\title {Isoperimetry of waists and local versus global asymptotic 
    convex geometries}
\author{Roman Vershynin}

\maketitle

\section{Introduction}
 
\begin{theorem}                 \label{two bodies}
  Assume that two symmetric convex bodies $K$ and $L$ in $\R^n$
  have sections of dimensions at least $k$ and $n-ck$ respectively 
  whose diameters are bounded by $1$.
  Then for a random orthogonal operator $U \in \OO(n)$
  the body $K \cap UL$ has diameter bounded by $C^{n/k}$
  with probability at least $1 - e^{-n}$.
\end{theorem}

\noindent 
Here and thereafter $C, C_1, c, c_1, \ldots$ denote positive absolute constants.

\qquad

{\bf 1.} The main point of Theorem \ref{two bodies} is that the 
existence implies randomness in it.
Namely, we do not assume that the two sections are random; 
their {\em existence} suffices. Yet the conclusion holds for a 
{\em random} rotation $U$.  
This seems to be a new phenomenon in the asymptotic convex geometry. 
It is further manifested by Corollary \ref{sections}.

{\bf 2.} Theorem \ref{two bodies} is a local to global
implication in the asymptotic convex geometry, see \cite{MS}. 
The {\em local} information about $K$ and $L$
(the existence of bounded sections) implies the {\em global} information
(bounded intersections of the whole bodies). 
This is further illustrated in Corollary \ref{global vr}.

{\bf 3.} The exponential bound $C^{n/k}$ in Theorem \ref{two bodies} 
can be improved to a polynomial bound, say $C(n/k)^2$, at the 
cost of decreasing the probability from $1-e^{-n}$ to $1-e^{-k}$.
This will be proved in the Appendix by Mark Rudelson and the author.

\qquad

We will first discuss two applications of Theorem \ref{two bodies} 
and then turn to the method used in its proof, 
which is rather general and whose main ingredient
is the recent ``isoperimetry of waists'' due to M.Gromov. 

The first immediate consequence of Theorem \ref{two bodies} 
is a striking statement ``existence implies randomness'' 
about the diameters of sections of symmetric convex bodies $K$: 
\begin{quote}
  {\em ``If $K$ has a nicely bounded section, \\
  then most sections of $K$ are nicely bounded''}
\end{quote}
(with a certain loss of the diameter as well as of the dimension).
This phenomenon was discovered by A.Giannopoulos, V.Milman 
and A.Tsolomitis in their forthcoming work \cite{GMT}, 
and independently by the author a few weeks later. 
Precisely, with $L = \R^{ck}$ or $K = \R^{n-k}$
in Theorem \ref{two bodies} one immediately obtains

\begin{corollary}[Propagation of boundedness of sections]          \label{sections}
  Let $K$ be a convex symmetric body in $\R^n$ and $k$ be
  a positive integer.

  (i) If there {\bf exists} a section of $K$ of dimension $k$ 
  whose diameter is bounded by $1$, then
  a {\bf random} section of $K$ of dimension $ck$
  has diameter bounded by $C^{n/k}$
  with probability at least $1 - e^{-cn}$.

  (ii) If there {\bf exists} a section of $K$ of dimension $n-ck$ 
  whose diameter is bounded by $1$, then
  a {\bf random} section of $K$ of dimension $n-k$
  has diameter bounded by $C^{n/k}$
  with probability at least $1 - e^{-cn}$.
\end{corollary}

\noindent The randomness here is with respect to the Haar measure on 
the Grassmanian $G_{n,m}$. 

The forthcoming paper \cite{GMT} offers a more direct approach to 
this corollary as well as better bounds on the diameter (note also 
that the version of Theorem \ref{two bodies} in the Appendix 
gives polynomial bounds). 

\setlength{\unitlength}{1pt}
\begin{center} 
\begin{picture}(140,140)(-20,-20)
\put(0,0) {\line(0,1){110}} 
\put(0,0) {\line(1,0){110}} 
\multiput(0,0)(20,20){5} {\line(1,1){10}}
\put(0,0) {\line(4,1){80}}
\put(80,20) {\line(1,4){20}}
\put(100,-5) {\line(0,1){10}} \put(96,-14){$n$}
\put(-5,100) {\line(1,0){10}} \put(-14,97){$n$}
\put(20,-14){{\small Existent}}
\put(-45,40){{\small Random}}
\end{picture} \\ \nopagebreak
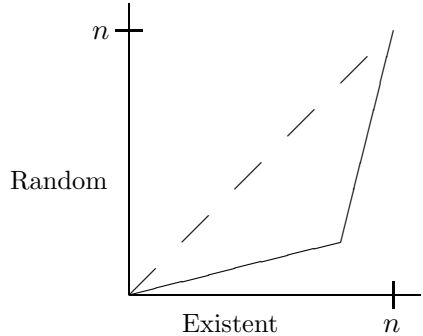
Dimensions of the existent and of the random sections
in Corollary \ref{sections}
\end{center}

\qquad

The second application is that Theorem \ref{two bodies} 
can turn various local results in the asymptotic convex geometry 
into global statements.
Let us show this on the example of the volume ratio theorem, 
one of the important ``local'' results in the field. 
A convex set $K$ in $\R^n$ has the volume ratio $A$ with respect
to the unit Euclidean ball $D$ if $D \subseteq K$ and 
$(|K|/|D|)^{1/n} = A$.

\begin{corollary}[Global volume ratio theorem]             \label{global vr}
  Assume that a convex set $K$ in $\R^n$ has volume ratio $A$
  with respect to the unit Euclidean ball. 
  Assume that a convex symmetric set $L$ in $\R^n$
  has a section of dimension $k$ whose diameter is bounded by $1$. 
  Then for a random orthogonal operator $U \in \OO(n)$
  the body $K \cap UL$ has diameter bounded by $(2A)^{Cn/k}$
  with probability at least $1 - e^{-n}$.
\end{corollary}

For $L = \R^{n-k}$, Corollary \ref{global vr}
is the classical volume ratio theorem due to S.Szarek and 
N.Tomczak-Jaegermann (see e.g. \cite{P}); 
the best constant in this case is known 
to be $C = 1$ (with $4\pi$ replacing the factor of $2$). 

\qquad

\proof By Rogers-Shephard \cite{RS}, the volume of 
$K' = K - K$ is 
$|K'| \le \binom{2n}{n} |K| \le 4^n |K|$.
Then $K'$ is symmetric and its volume ratio with respect to 
the Euclidean ball is at most $4A$. By the volume ratio theorem 
(see e.g. \cite{P}), $K'$ has a section of dimension at least
$n-ck$ whose diameter is bounded by $M = (4\pi A)^{n/ck}$.
The proof is finished by applying Theorem \ref{two bodies} to 
$M^{-1} K'$ and $L$.
\endproof  

\qquad

Our approach to the proof of Theorem \ref{two bodies}
is based on a recent isoperimetric theorem of M.Gromov \cite{G}, 
his ``isoperimetry of waists'' on the unit Euclidean sphere $S^{n-1}$:
\begin{quote}
  If $f : S^{k-1} \to S^{n-1}$ is an odd and continuous map, 
  then the $(n-1)$-volume of any $\e$-neighborhood of $f(S^{k-1})$
  in the geodesic distance on the sphere is minimized when 
  $f$ is the canonical embedding, i.e. when the ``waist'' $f(S^{k-1})$
  is an equatorial sphere.
\end{quote}
This is proved in \cite{G} for certain $k$ and $n$ and it remains
an open problem for the rest of $k,n$; see next section.

The isoperimetry of waists can be effectively used in the asymptotic 
convex geometry. Suppose we know that for a symmetric convex body 
$K$ in $\R^n$ {\em there exists an orthogonal projection $PK$ 
that contains the unit Euclidean ball}.
Without loss of generality, let $S^{k-1}$ be the sphere of that ball.
One can find an odd and continuous lifting $g : S^{k-1} \to K$ of the projection 
$P$ and contract it to the sphere by defining $f(x) = g(x)/|g(x)|$.
Then $f : S^{k-1} \to S^{n-1}$ satisfies the assumptions of Gromov's 
isoperimetry and, moreover, the waist $f(S^{k-1})$ lies in $K$.
Then the isoperimetry of waists gives a {\em computable lower bound
on the $(n-1)$-volume of any $\e$-neighborhood of $K$ on the sphere $S^{n-1}$}.
This bound is sharp; it reduces to an equality if the projection 
$PK$ coincides with the section $K \cap P\R^{n}$. The exact statement
is Proposition \ref{projection}.
 
This argument is the main step in the proof of (the dual form of) 
Theorem~\ref{two bodies}. The assumptions are that both $K$ and $L$ 
have orthogonal projections that contain unit Euclidean balls. 
The reasoning above based on the isoperimetry of waists implies
that the appropriate neighborhoods $K_{\e_1}$ of $K$ and 
$L_{\e_2}$ of $L$ have large $(n-1)$-volumes on the unit sphere.
Then a standard $\e$-net argument shows that the Minkowski sum 
$K_{\e_1} + U L_{\e_2}$ contains the unit Euclidean ball with large
probability (see Lemma \ref{core}). 
If $\e_1 + \e_2$ is a small number, then $K + UL$ must 
contain some nontrivial Euclidean ball, too. This is (the dual form of) 
the conclusion of Theorem~\ref{two bodies}.

\qquad

\noindent ACKNOWLEDGEMENTS.
The author is grateful to Apostolos Giannopoulos,
Vitali Milman and Gideon Schechtman for important discussions.

\section{Gromov's isoperimetry of waists}

The normalized Lebesgue measure on the unit Euclidean sphere 
$S^{n}$ will be denoted by $\s_n$. For a subset $A \subset S^n$ and 
a number $\theta > 0$, by $A_\theta$ we denote the $\theta$-neighborhood
of $A$ in the geodesic distance $d$, i.e. 
$A_\theta = \{ y \in S^n :\; d(x,y) \le \theta \text{ for some } x \in A \}$.
A map $f$ is called odd is $f(-x) = -f(x)$ for all $x$.
The following isoperimetry is proved in \cite{G} 6.3.B.

\begin{theorem}[Gromov's isoperimetry of waists]
  Let $n$ be odd and $l = 2^k - 1$ for some integer $k$.
  Let $f: S^{n-l} \to S^n$ be an odd continuous function.
  Then for all $0 < \theta < \pi/2$ 
  $$
  \s_{n}((f(S^{n-l}))_\theta) \ge \s_n((S^{n-l})_\theta).
  $$
\end{theorem}
  
Conjecturally, this theorem should hold for all $n,l$. 
We will actually need this for all $n,l$, 
and in the absence of such result we will deduce a
relaxed version of Gromov's theorem for all $n,l$, 
Corollary \ref{general Gromov} below. 
This is done naturally by embedding into a higher dimensional sphere.

\begin{lemma}                   \label{higher sphere}
  Let $A \subset S^n$ be a symmetric measurable set
  and $m \ge n$ be a positive integer. 
  Then for all $0 < \theta < \pi/2$
  $$
  \s_n(A_\theta)  \ge  \s_m(A_\theta),
  $$
  where in the right side we look at a set $A$ as a subset of $S^m$
  via the canonical embedding $S^n \subset S^m$.
\end{lemma}

\proof
Fix an $x \in S^{m}$ and let $x_1$ be its spherical projection 
onto $S^n$, i.e. $x_1 = P_n x / |P_n x|$, where $P_n$ denotes
the orthogonal projection in $\R^{m+1}$ onto $\R^{n+1}$.

\vspace{0.2cm} 

CLAIM: $d(x_1, A) \le d(x, A)$. 

\vspace{0.2cm}

\noindent 
To prove the claim, since $A$ is symmetric it is enough to check that 
$$
\text{$d(x_1,a) \le d(x,a)$ 
      for all $a \in A$ such that $d(x,a) \le \pi/2$.}
$$ 
Since $0 \le d(x,a) \le \pi/2$ and $\< x, a\> = \cos d(x,a)$, 
we have 
$$
0 \le \< x,a\> \le 1.
$$
Since $a \in S^n$, we have $P_n a = a$; thus 
$$
\< x_1,a\> = \< P_n x / |P_n x|, a\> 
= \frac{1}{|P_n x|} \< x, a\> 
\ge \< x,a\> .
$$
In particular, $0 \le \< x_1,a\> \le 1$. 
Since the function $\cos^{-1} : [0,1] \to [0,\pi/2]$ is decreasing, 
$$
d(x_1,a) = \cos^{-1} \< x_1,a\> 
\le \< x,a\> = d(x,a).
$$
This proves the Claim.

\vspace{0.2cm}

Now we can finish the proof of the lemma as follows: 
\begin{align*}
\s_m(A_\theta) 
  &=   \s_m(x \in S^m \;:\; d(x,A) \le \theta) \\
  &\le \s_m(x \in S^m \;:\; d(x_1,A) \le \theta) \\
  &= \s_n(x_1 \in S^n \;:\; d(x_1,A) \le \theta)
   = \s_n(A_\theta)
\end{align*}
where the last line is obtained by representing a uniformly 
distributed vector $x \in S^m$ as $x = \g x_1 + \sqrt{1-\g} x_2$, 
where $x_1 \in S^n$ and $x_2 \in S^{m-n}$ are uniformly distributed, 
$\g$ is an appropriate random variable and the three random variables 
$x_1, x_2, \g$ are jointly independent. 
\endproof

\begin{corollary}[General (relaxed) isoperimetry of waists]   \label{general Gromov}
  Let $l < n$ are positive integers.
  Let $f: S^{n-l} \to S^n$ be a n odd continuous function.
  Then for all $0 < \theta < \pi/2$ 
  $$
  \s_{n}((f(S^{n-l}))_\theta) \ge \s_{n+l+1}((S^{n-l-1})_\theta).
  $$
\end{corollary}
  
\proof

CASE 1: $n-l$ is even. 

Let $k$ be the minimal integer such that $2^k - 1 \ge l$. 
Then $m := (n-l) + (2^k-1)$ is odd. 
Moreover, since $2^{k-1} < l$, we have $2^k \le 2(l+1)$, so 
$m < n - l + 2(l+1) - 1  \le  n + l + 1$. Hence 
$$
n \le m \le n+l.
$$
Then Gromov's theorem can be applied to functions 
from $S^{n-l} \to S^{m}$, in particular to 
$f : S^{n-l} \to S^n \to S^m$
where the second map is the canonical embedding. 
Then using Lemma \ref{higher sphere}, Gromov's theorem and 
Lemma \ref{higher sphere} again, we have
$$
\s_{n}(f(S^{n-l})_\theta)
\ge \s_{m}(f(S^{n-l})_\theta)
\ge \s_{m}(S^{n-l}_\theta)
\ge \s_{n+l}(S^{n-l}_\theta).
$$

CASE 2: $n-l$ is odd.

Apply Case 1 to the function $g: S^{n-l-1} \to S^{n-l} \to S^n$
where the first map is the canonical embedding and the second map 
is $f$. We have
$$
\s_{n}(f(S^{n-l})_\theta) 
\ge \s_{n}(g(S^{n-l-1})_\theta) 
\ge \s_{n+l+1}(S^{n-l-1}_\theta).
$$
Therefore for all $l < n$ we have
$$
\s_{n}(f(S^{n-l})_\theta) \ge \s_{n+l+1}(S^{n-l-1}_\theta)
$$
(here we used Lemma \ref{higher sphere} again to step one dimension 
up in Case 1).
\endproof

To simplify the use of Corollary \ref{general Gromov}, 
we will denote:
$$
\s_{n,k}(\theta) = \s_n((S^k)_\theta), \ \ \ 
\s^\Lip_{n,k}(\theta) = \inf_f \s_n((f(S^k))_\theta), 
$$
where the infimum is over all symmetric continuous functions 
$f: S^k \to S^n$.

\begin{general Gromov'} 
  Let $k < n$ be positive integers. Then for $0 < \theta < \pi/2$
  \begin{equation}                  \label{gen Gromov}
  \s^\Lip_{n,k}(\theta)  \ge  \s_{2n-k+1, k-1}(\theta).
  \end{equation}
\end{general Gromov'}

\remark If Gromov's theorem is true for all $n,l$, then 
Corollary \ref{general Gromov}$'$ improves to 
$$
\s^\Lip_{n,k}(\theta) \ge \s_{n,k}(\theta).
$$

\qquad

The right hand side of \eqref{gen Gromov} is a computable quantity.
Sharp asymptotic estimates on $\s_{n,k}(\theta)$ were found by 
S.Artstein \cite{A}. For our present purpose, we will be satisfied 
with less precise estimates, which reduce to computations
on Gaussians and whose prove we include for completeness.

\begin{lemma}                       \label{cap large}
  Let $1 < k \le n$ be integers and let $0 < \e < 1/2$. Then 
  $$
  (c\e)^{2k}  
  \le  \s_{n-1,n-k-1}(\sin^{-1}\sqrt{\frac{\e^2 k}{n}})  
  \le  (C\e)^{k/2}.
  $$
  Consequently,
  $$
  1 - (C\e)^{k/2}  
  \le  \s_{n-1, k-1}(\sin^{-1}\sqrt{1-\frac{\e^2 k}{n}})  
  \le  1 - (c\e)^k.
  $$
\end{lemma}
For the proof, we quote two known facts about the 
canonical real Gaussian vector.

\begin{fact}                        \label{gaussian}
  Let $g_1, g_2, \ldots$ be a sequence of i.i.d. normalized 
  Gaussian random variables. Then 

  (i) For every $M \ge 2$ one has 
      $$
      \P \{ g_1^2 + \cdots + g_k^2 > M^2 k \} 
      \le 2e^{-c M^2 k};
      $$

  (ii) for every $\e> 0$, we have 
      $$
      (c \e)^k \le \P \{ g_1^2 + \cdots + g_k^2 \le \e^2 k \} 
      \le (C \e)^k.
      $$
\end{fact}

\noindent{\bf Proof of Lemma \ref{cap large}. }
(i) By the rotation invariance of the Gaussian density, 
$$
\s := \s_{n-1,n-k-1}(\sin^{-1}\sqrt{\frac{\e^2 k}{n}}) 
= \P \{ (g_1^2 + \cdots + g_k^2) 
         \le \frac{\e^2 k}{n} (g_1^2 + \cdots + g_n^2) \}.
$$
If we write $\frac{\e^2 k}{n} = \frac{\e^4 k}{\e^2 n}$ 
then by Fact \ref{gaussian} (ii) we have 
$$
\s \ge \P \{ g_1^2 + \cdots + g_k^2 \le \e^4 k \}
     - \P \{ g_1^2 + \cdots + g_n^2 \le \e^2 n \}
\ge (c\e^2)^k - (C\e)^n \ge (c_1\e)^{2k}
$$
since $1 < k < n/4$.

To prove the reverse inequality, let $M \ge 2$ and write 
$\frac{\e^2 k}{n} = \frac{M^2 \e^2 k}{M^2 n}$. 
Then by Fact \ref{gaussian} (i) and (ii) we have 
\begin{equation}                    \label{sigma via M}
\s \le \P \{ g_1^2 + \cdots + g_k^2 \le M^2 \e^2 k \}
     - \P \{ g_1^2 + \cdots + g_n^2 > M^2 n \}
\le (CM\e)^k + 2 e^{-c M^2 n}.
\end{equation}
If $e^{-4cn} \le (2C\e)^k$ then letting $M=2$ in \eqref{sigma via M}
we obtain $\s \le 3(C\e)^k$, as required.
Thus we can assume that 
$e^{-4cn} > (2C\e)^k > \e^{4k}$, hence $\log(1/\e) \ge (cn/k)$.
Let $M = 2\sqrt{(k/cn) \log(1/\e)}$. Note that 
$2 < M \le C_1 \sqrt{\log(1/\e)}$. 
With this $M$ in \eqref{sigma via M} we obtain
$$
\s \le (C_2 \e \sqrt{\log(1/\e)})^k + 2\e^{-4k}
\le (C_3\e)^{k/2}.
$$
This proves the first part of the Lemma. 
The second part follows from the equation 
$$
\s_{n-1,n-k-1}(\sin^{-1}\alpha) 
+ \s_{n-1, k-1}(\sin^{-1}\sqrt{1-\alpha^2}) 
= 1, 
\ \ \ 0 < \a < 1,
$$
and the first part. 
\endproof

When the general Gromov's theorem (Corollary \ref{general Gromov}$'$)
is combined with Lemma \ref{cap large}, we obtain explicit estimates
for $\s^\Lip_{n,k}(\theta)$:

\begin{corollary}               \label{large lip}
  Let $1 < k \le n$ be integers and let $0 < \e < 1/2$. Then 
  
  (i)  $\s^\Lip_{n-1,n-k-1}(\sin^{-1}\sqrt{\frac{\e^2 k}{n}})  
       \ge  (c\e)^{8k}$;

  (ii) $\s^\Lip_{n-1, k-1}(\sin^{-1}\sqrt{1-\frac{\e^2 k}{n}})  
       \ge  1 - (C\e)^{k/4}$.
\end{corollary}

\proof
Let $\a = k/n$.

(i) By Corollary \ref{general Gromov}$'$, 
\begin{equation}                    \label{large lip 1}
\s^\Lip_{n-1,(1-\a)n-1}(\sin^{-1}\sqrt{\e^2\a})  
       \ge \s_{n+\a n,(1-\a)n-2}(\sin^{-1}\sqrt{\e^2\a}).
\end{equation}
To apply Lemma \ref{cap large}, write the right hand side of 
\eqref{large lip 1} for suitable $m$ and $\b$ as
\begin{equation}                    \label{large lip 2}
\s_{m-1, (1-\b)m-1} (\sin^{-1} \sqrt{(\e^2\a/\b) \cdot \b})
\ge (c \sqrt{\e^2\a/\b})^{2\b m}.
\end{equation}
The numbers $m$ and $\b$ are, of course, determined by 
$m-1 = n + \a n$ and $(1-\b)m - 1 = (1-\a)n - 2$. Hence
$\b = (2\a n + 2)/(n+ \a n + 1)$, so that $\a < \b < 3\a$.
Then we can continue \eqref{large lip 2} as 
$$
\ge (c_1\e)^{4(\a n + 1)} \ge (c_1 \e)^{8k}.
$$
This completes part (i).

(ii)  By Corollary \ref{general Gromov}$'$, 
\begin{equation}                    \label{large lip 3}
\s^\Lip_{n-1, \a n-1}(\sin^{-1}\sqrt{1-\e^2\a})  
       \ge \s_{2n-\a n, \a n-2}(\sin^{-1}\sqrt{1-\e^2\a}).
\end{equation}
To apply Lemma \ref{cap large}, write the right hand side of 
\eqref{large lip 3} for suitable $m$ and $\b$ as
\begin{equation}                    \label{large lip 4}
\s_{m-1, \b m-1} (\sin^{-1} \sqrt{1 - (\e^2\a/\b) \cdot \b})
\ge 1 - (C \sqrt{\e^2\a/\b})^{\b m} - e^{-10m}.
\end{equation}
The numbers $m$ and $\b$ are, of course, determined by 
$m-1 = 2n - \a n$ and $\b m - 1 = \a n - 2$. Hence
$\b = (\a n - 1)/(2n - \a n + 1)$, so that $\b \ge \a/2$.
Then we can continue \eqref{large lip 4} as 
$$
\ge 1 - (C_1 \e)^{(\a n - 1)/2} 
\ge 1 - (C_1 \e)^{k/4}.
$$
This completes part (ii).
\endproof

\section{Waists generated by projections of convex bodies}
The following observation connects the isoperimetry of waists to 
convex geometry.

For simplicity, given a set $A \in \R^n$ we write 
$\s_{n-1} (A)$ for $\s_{n-1} (A \cap S^{n-1})$, if measurable. 
The unit Euclidean ball in $\R^n$ is denoted by $D$. 
Minkowski sum in $\R^n$ is defined as 
$A + B = \{ a+b :\; a \in A, \; b \in B \}$.

\begin{proposition}                 \label{projection}
  Let $K$ be a convex symmetric set in $\R^n$. 
  Assume there is an orthogonal projection $P$, 
  $\rank P = k$, such that $PK \supseteq PD$.
  Then for all $0 < \e < 1$
  \begin{equation}              \label{projection est}
  \s_{n-1} (K + \e D)  
  \ge \s^\Lip_{n-1,k-1}(\sin^{-1}\e).
  \end{equation}
\end{proposition}

\remark 
The power of this fact is that the right side of \eqref{projection est} is 
easily estimated via Gromov's theorem (Corollary \ref{large lip}).

\qquad

\proof
We can assume that the range of $P$ is $\R^k$, so 
$PK \supseteq S^{k-1}$.
There exists an odd continuous lifting $g : S^{k-1} \to K$ of the 
projection $P$. Define 
$$
f : S^{k-1} \to S^{n-1}, \ \ \ f(x) = g(x) / |g(x)|.
$$
The function $f$ is odd and it is continuous because
$$
|g(x)| \ge |Pg(x)| = |x| = 1
\ \ \ \text{for all $x$.}
$$
Since also $g(x) \in K$, we have $f(x) \in K$, thus 
$$
f(S^{k-1}) \subseteq K \cap S^{n-1}.
$$
By making a simple planar drawing, one sees that 
for every $y \in S^{n-1}$
\begin{equation}                \label{interval}
[-y,y] + \e D  \supseteq \{y\}_{\sin^{-1}\e}
\end{equation}
Running $y$ over $f(S^{k-1})$, we obtain 
\begin{align*}
K + \e D  
  &\supseteq  f(S^{k-1}) + \e D \\
  &= \bigcup_{y \in f(S^{k-1})} \Big( [-y,y] + \e D \Big)
       \ \ \text{by the symmetry of $f$} \\
  &= f(S^{k-1})_{\sin^{-1}\e}
       \ \ \text{by \eqref{interval}.} 
\end{align*} 
Intersecting both sides with $S^{n-1}$ and taking the 
measure completes the proof. 
\endproof

Proposition \ref{projection} will in particular be used
to estimate the covering number of $K + \e D$.

Given two convex sets $L$ and $K$, the covering number $N(L,K)$ 
is the minimal number of translates of $K$ needed to cover $L$.
By a simple and known volumetric argument, 
$N(L,K) \le \frac{|L+K|}{|K|}$. 

\begin{lemma}                   \label{entropy}
  For every convex symmetric set $K$, 
  $$
  N(D,K)  \le  2^n / \s_{n-1}(K).
  $$
\end{lemma}

\proof
\begin{equation}                \label{entropy small}
N(D,K)  \le N(D, K \cap D)
\le  \frac{|D + (K \cap D)|}{|K \cap D|}
\le  \frac{|2 D|}{|K \cap D|}.
\end{equation}
Next, 
\begin{equation}                \label{surface small}
\s_{n-1}(K) = \s_{n-1}(K \cap D)
\le \frac{|K \cap D|}{|D|},
\end{equation}
which folows from a standard argument that transfers the 
surface measure on $S^{n-1}$ to the volume in $D$
(a set $A \subseteq S^{n-1}$ generates the cone $\cup_{0<t<1} tA$,
which occupies the same portion of the volume in $D$ 
as $\s_{n-1}(A)$).  

Then \eqref{entropy small} and \eqref{surface small} complete 
the proof.
\endproof

\section{Proof of Theorem \ref{two bodies}}

By duality, Theorem \ref{two bodies} can equivalently be stated as follows. 
There exist an absolute constant $a \in (0,1)$ such that the following holds. 
Assume that there exist orthogonal projections $P$ and $Q$ with 
$\rank P = k$ and $\rank Q = n-ak$, and such that 
\begin{equation}                        \label{PQ}
  PK \supseteq PD, \ \ \ QL \supseteq QD.
\end{equation}
Then for $U$ as in the theorem, we claim that
\begin{equation}                        \label{WTS}
  K + UL \supseteq C^{n/k} D.
\end{equation}

The idea is as follows. Let $\d_K, \d_L > 0$ be parameters.
By Gromov's theorem and Lemma \ref{entropy}, we will be able 
to estimate
\begin{equation}                        \label{sigma N}
  1-\s := \s_{n-1} (K + \d_K D)
  \ \ \ \text{and} \ \ \ 
  N := N(2D, L + \d_L D).
\end{equation}

\begin{lemma}                           \label{core}
  Let $K$ and $L$ be convex bodies in $\R^n$ such that 
  \eqref{sigma N} holds and $\d_K + \d_L < 1$. 
  Then for a random orthogonal operator $U \in \OO(n)$
  \begin{equation}                      \label{core eq}
    (1 - \d_K - \d_L) D  \subseteq K + UL
  \end{equation}
  with probability at least $1 - N\s$.
\end{lemma}

\proof
By a standard argument, the sphere $S^{n-1}$ of $D$ can be covered
by $N$ translates of the body $L + \d_L D$ by vectors from $S^{n-1}$.
Hence there exists a subset $\NN \subset S^{n-1}$ such that 
\begin{equation}                        \label{NN}
|\NN| = N 
\ \ \ \text{and} \ \ \ 
D \subseteq \NN + L + \d_L D.
\end{equation}
Since for every $z \in S^{n-1}$, its image $Uz$ under a random rotation 
$U \in \OO(n)$ is uniformly distributed on the sphere, we have
for any fixed $z \in \NN$:
$$
\P \{ U \in \OO(n) \;:\; Uz \in K + \d_K D \}
= \s_{n-1}(K + \d_K D) 
= 1 - \s.
$$
Thus
$$
\P \{ U \in \OO(n) \;:\; U\NN \subseteq K + \d_K D \}
\ge 1 - N \s.
$$
Fix any $U$ in this set and apply it to the inclusion in \eqref{NN}:
$$
D
\subseteq U\NN + UL + \d_L D
\subseteq K + \d_K D + UL + \d_L D.
$$
Since $\d_K + \d_L < 1$, this inclusion implies \eqref{core eq}.
\endproof

\qquad

\noindent {\bf Proof of Theorem \ref{two bodies}. }
We can clearly assume that $0 < a < 1/33$ and that $ak \ge 1$.
Let
$$
\e_K > 0, \ \ \ \d_K = \sqrt{1 - \frac{\e_K^2 k}{n}}.
$$
By Proposition \ref{projection} and Corollary \ref{large lip} (ii), 
$$
\s_{n-1}(K + \d_K D)
\ge  \s^\Lip_{n-1,k-1}(\sin^{-1} \d_K)\\
\ge 1 - (C \e_K)^{k/4}
\ge 1 - 2 e^{-10n}
$$
if one chooses the value of $\e_K$ as
$$
\e_K = \exp(-C_1 n/k), 
$$
where $C_1 > 0$ is a sufficiently large absolute constant.
Similarly, let 
$$
\e_L > 0, \ \ \ \d_L = \sqrt{\frac{\e_L^2 ak}{n}}.
$$
By Proposition \ref{projection} and Corollary \ref{large lip} (i), 
$$
\s_{n-1}(L + \d_L D)
\ge \s^\Lip_{n-1, n-ak-1}(\sin^{-1} \d_L) \\
\ge (c\e_L)^{8k}
\ge \frac{1}{2} e^{-n/2}
$$
if one chooses the value of $\e_L$ as
$$
\e_L = \exp(-c_2 n/ak)
$$
where $c_2 > 0$ is a sufficiently small absolute constant.
By Lemma \ref{entropy}, 
$$
N(2D, 2L + 2\d_L D) 
= N(D, L + \d_L D) 
\le 2^n / \frac{1}{2} e^{-n/2}
\le 2 e^{1.2n}.
$$ 
By Lemma \ref{core}, if $\d_K + 2\d_L < 1$ then the desired inclusion
\begin{equation}                        \label{desired incl}
(1 - \d_K - 2\d_L) D  \subseteq K + 2UL
\end{equation}
holds with probability at least
$$
1 - 2e^{1.2n} \cdot 2e^{-10n}
\ge 1 - e^{-n}.
$$
So it only remains to bound below
$$
\d_K + 2\d_L 
= \sqrt{1 - \exp(-2 C_1 n/k) (k/n)} + 2\sqrt{\exp(-2 c_2 n/ak) (ak/n)}.
$$
This can be quickly done using the inequalities 
$\sqrt{1-x} \le 1 - x/2$ and $x e^{-C/x} \ge e^{-2C/x}$
valid for all $0 < x < 1$ and for a sufficiently large absolute 
constant $C$. We thus have 
$$
\d_K + 2\d_L 
\le 1 - \frac{1}{2} \exp(-C'_1 n/k) + 2\exp(-c'_2 n/ak)
< 1 - \frac{1}{4} \exp(-Cn/k)
$$
if $a$ is chosen a sufficiently small absolute constant.
This together with \eqref{desired incl} completes the proof.
\endproof

\end{document}